\documentclass{article}
\usepackage[utf8]{inputenc}

\usepackage{amsmath, amsthm, amssymb, mathrsfs, thmtools, enumerate}
\usepackage[colorlinks = true, linkcolor = black, urlcolor = blue, citecolor = black, pdfborder={0 0 0}]{hyperref} 
\usepackage{tikz} 

\usepackage{xcolor}
\definecolor{junglegreen}{RGB}{65,150, 00}
\definecolor{savanahblue}{RGB}{0,65,150}
\definecolor{saharahred}{RGB}{150,20,40}
\allowdisplaybreaks[4]
\declaretheoremstyle[
  spaceabove=1em, spacebelow=1em,
  headfont= \bfseries,
  notefont=\mdseries, notebraces={(}{)},
  bodyfont=\normalfont,
  postheadspace=1em,
  qed= \tiny $\blacksquare$
]{example}

\declaretheoremstyle[
  spaceabove=\topsep, spacebelow=\topsep,
  headfont= \bfseries,
  notefont=\mdseries, notebraces={(}{)},
  bodyfont=\normalfont,
  postheadspace=1em,
  qed= $\blacktriangledown$
]{remark}

\declaretheoremstyle[
  spaceabove=\topsep, spacebelow=\topsep,
  headfont= \bfseries,
  notefont=\mdseries, notebraces={(}{)},
  bodyfont=\it{\normalfont},
  postheadspace=1em,
  qed= \qedsymbol
]{obvious}

\theoremstyle{plain}
\newtheorem{Thm}{Theorem}[section]

\newtheorem{Cor}[Thm]{Corollary}
\newtheorem{Lem}[Thm]{Lemma}

\theoremstyle{definition}

\declaretheorem[style=example, name = Example, sharenumber = Thm]{Exam}

\numberwithin{equation}{section}
\numberwithin{table}{section}

\newcommand{\N}{\mathbb{N}}

\newcommand{\Q}{\mathbb{Q}}
\newcommand{\QG}{\Q[G]}
\newcommand{\QH}{\Q[H]}

\renewcommand{\S}{\mathfrak{S}}
\newcommand{\T}{\mathfrak{T}}
\newcommand{\z}{\mathcal{Z}}
\renewcommand{\H}{\mathcal{H}}
\renewcommand{\L}{\mathcal{L}}
\newcommand{\K}{\mathcal{K}}
\newcommand{\D}{\mathcal{D}}
\newcommand{\dsum}{\displaystyle\sum}

\DeclareMathOperator{\Aut}{Aut}
\DeclareMathOperator{\Span}{Span}

\newcommand{\thmref}[1]{Theorem \ref{#1}}
\newcommand{\corref}[1]{Corollary \ref{#1}}
\newcommand{\lemref}[1]{Lemma \ref{#1}}

\newcommand{\tableref}[1]{Table \ref{#1}}
\newcommand{\secref}[1]{Section \ref{#1}}

\begin{document}

\title{Counting Schur Rings over Cyclic Groups of Semiprime Order}
\author{Joseph Keller, Andrew Misseldine\footnote{Andrew Misseldine, Southern Utah University, andrewmisseldine@suu.edu, (telephone) 435-865-8228, (fax) 435-865-8666}, Max Sullivan}
\date{\today}

\maketitle


\begin{abstract}
In this paper, we continue the enumeration of Schur rings over cyclic groups. Cyclic groups of semiprime order $pq$, where $p$ and $q$ are distinct primes, are considered. Additionally, cyclic groups of order $4p$ are considered. 



\end{abstract}

\section{Introduction}
Let $G$ be a finite group, and let $\QG$ denote the rational group algebra. Let $\L(G)$ denote the lattice of subgroups of $G$. 
For any subset $C\subseteq G$, we may associate to it an element of the group algebra, namely $\sum_{g\in C} g\in \QG$. Such an element is called a \emph{simple quantity}. When the context is clear, we will likewise denote the simple quantity associated to $C$ as $C$ itself. Define 
$C^* := \{x^{-1} \mid x\in C\}$ for all 
$C\subseteq G$. Let $\{C_1, C_2, \ldots, C_r\}$ be a partition of $G$, and let $\S$ be the subspace of $\QG$ spanned by $\{C_1, C_2, \ldots, C_r\}$. We say that $\S$ is a \emph{Schur ring}\label{def:schurring} over $G$ if 
\begin{enumerate}
\item $C_1 = \{1\}$; 
\item For each $i$, there is a $j$ such that $C_i^* = C_j$;
\item For each $i$ and $j$, $C_iC_j = \dsum_{k=1}^r \lambda_{ijk}C_k$, for $\lambda_{ijk}\in \N$.
\end{enumerate}
The sets $C_1, C_2, \ldots, C_r$ are called the \emph{$\S$-classes} (or \emph{primitive sets} of $\S$). Note that a Schur ring is uniquely determined by its associated partition of $G$. We will denote this partition as $\D(\S)$. 

Schur rings over cyclic groups have been of great interest for the last few decades because of their connection to algebraic graph theory and association schemes (see \cite{Muzychuk09}). In \cite{CountingII}, the second author develops a general technique of enumerating cyclic Schur rings that we will employ in the case of $n=pq, 4p$, where $p$ and $q$ are distinct primes. This paper is then a continuation of work begun by the second author in \cite{CountingII, MePhD, Counting}.

Let $\z_n =\langle z\rangle$ denote the cyclic group of order $n$. Let $\Omega(n)$ denote the number of Schur rings over $\z_n$. In this paper, we provide formulas for $\Omega(pq)$ and $\Omega(4p)$. We present first the semiprime case $pq$:

\begin{Thm}\label{thm:main} Let $p$ and $q$ be distinct primes such that $p=\prod_{i=1}^n r_i^{k_i} +1$ and $q=\prod_{i=1}^n r_i^{\ell_i} +1$, where $\{r_1, r_2,\ldots, r_n\}$ is a list of distinct primes. Then 
\[\Omega(pq) = \prod_{i=1}^n \sum_{j=0}^{\min(k_i,\ell_i)} \phi(r_i^{j})(k_i-j+1)(\ell_i-j+1) + 2\prod_{i=1}^n (k_i+1)(\ell_i+1) + 1,\] where $\phi$ denotes Euler's totient function.
\end{Thm}

Be aware that in the above decompositions of $p$ and $q$, the exponents $k_i$ and $\ell_i$ may possibly be zero, allowing for a common family of primes $\{r_1, r_2, \ldots, r_n\}$ between $p$ and $q$.

We list next some useful simplifications of \thmref{thm:main} when special conditions are placed on the primes $p$ and $q$. The proofs of the following corollaries are immediate from \thmref{thm:main}.

\begin{Cor}\label{cor:main} Let $p$ and $q$ be distinct primes such that $p=2^ka +1$ and $q=2^\ell b+1$, where $a$ and $b$ are both odd integers and $\gcd(a,b)=1$. Let $x$ and $y$ be the number of divisors of $p-1$ and $q-1$, respectively. Then 
\[\Omega(pq) = \left[3(k+1)(\ell+1) + \sum_{j=1}^{\min(k,\ell)} 2^i(k-j+1)(\ell-j+1)\right]\left(\dfrac{xy}{(k+1)(\ell+1)}\right) + 1.\]
\end{Cor}

\begin{Cor}\label{cor:2p} Let $p\neq 2$ be a prime, and let $x$ be the number of divisors of $p-1$. Then 
\[\Omega(2p) = 3x+1.\]
\end{Cor}

\corref{cor:main} is particular useful when $p$ is a \emph{Fermat prime}, that is, $p=2^k+1$. There are only five know Fermat primes: $3$, $5$, $17$, $257$, and $65537$. It is widely conjectured that these are the only Fermat primes. We illustrate the simplification of \corref{cor:main} for the Fermat primes $3$ and $5$. 

\begin{Cor}\label{cor:3p} Let $p\neq 3$ be a prime such that $p=2^ka+1$ where $a$ is odd, and let $x$ be the number of divisors of $p-1$. Then 
\[\Omega(3p) = \left(\dfrac{7k+6}{k+1}\right)x + 1\] When $p\equiv 3 \pmod 4$, $\Omega(3p) = \dfrac{13}{2}x+1$.
\end{Cor}

\begin{Cor}\label{cor:5p} Let $p\neq 5$ be a prime such that $p=2^ka+1$ where $a$ is odd, and let $x$ be the number of divisors of $p-1$. Then 
\[\Omega(5p) = \left(\dfrac{13k+7}{k+1}\right)x + 1\] When $p\equiv 3 \pmod 4$, $\Omega(5p) = 10x+1$.
\end{Cor}

We mention that \corref{cor:main} is also applicable when $p$ is a \emph{safe prime}, that is, $p=2r+1$, where $r$ is itself a prime\footnote{In this case, $r$ is necessarily a Sophie Germain prime.}. It is widely conjectured that there are infinitely many safe primes, the first few being:\footnote{We have intentionally omitted $5$ from the list of safe primes as it is the only safe prime which is Fermat. As a consequence, it is the only safe prime $p$ for which the number of divisors of $p-1$ is 3 instead of 4.} $7$, $11$, $23$, $47$, $59$, $83$, and $107$. Note that by \corref{cor:main}, if $p$ and $q$ are both safe primes, then $\Omega(pq) = 53$. Likewise, if $p$ is a safe prime, then $\Omega(2p)=13$, $\Omega(3p) = 27$ and $\Omega(5p) = 41$, by Corollaries \ref{cor:2p}, \ref{cor:3p}, and \ref{cor:5p}, respectively.

Using Corollaries \ref{cor:2p}, \ref{cor:3p}, and \ref{cor:5p} and the above discussion of safe primes, one can easily compute the number of Schur rings over $\z_{pq}$ for all semiprimes under 100. These are listed in \tableref{table:NumberSRings}. The one exception here is $n=91=7\cdot 13 = (2\cdot 3+1)(2^2\cdot 3+1)$. In this case, $\Omega(91)$ can be computed directly using \thmref{thm:main}.

\begin{table}[hbt]
\caption{Number of Schur Rings over $Z_{pq}$}
\begin{center}
\begin{tabular}{|c|c||c|c||c|c||c|c||c|c|}
\hline
\textbf{\emph{n}} & $\boldsymbol{\Omega(n)}$ & \textbf{\emph{n}} & $\boldsymbol{\Omega(n)}$ & \textbf{\emph{n}} & $\boldsymbol{\Omega(n)}$ & \textbf{\emph{n}} & $\boldsymbol{\Omega(n)}$ & \textbf{\emph{n}} & $\boldsymbol{\Omega(n)}$  \\\hline
6 & 7 & 26 & 19 & 46 & 13 & 65 & 67 & 86 & 25 \\ \hline
10 & 10 & 33 & 27 & 51 & 35 & 69 & 27 & 87 & 41 \\ \hline
14 & 13 & 34 & 16 & 55 & 41 & 74 & 28 & 91 & 97 \\ \hline
15 & 21 & 35 & 41 & 57 & 40 & 77 & 53 & 93 & 53 \\ \hline
21 & 27 & 38 & 19 & 58 & 19 & 82 & 25 & 94 & 13 \\ \hline
22 & 13 & 39 & 41 & 62 & 25 & 85 & 60 & 95 & 61 \\ \hline
\end{tabular}
\end{center}
\label{table:NumberSRings}
\end{table}

We next present the counting formula for $4p$:

\begin{Thm}\label{thm:main2} Let $p$ be an odd prime such that $p = 2^ka+1$, where a is an odd integer and $x$ the number of divisors of $p-1$. Then \[\Omega(4p)=\frac{15k+14}{k+1}x+3.\]  
\end{Thm}

\begin{table}[hbt]
\caption{Number of Schur Rings over $Z_{4p}$}
\begin{center}
\begin{tabular}{|c|c||c|c||c|c||c|c||}
\hline
\textbf{\emph{n}} & $\boldsymbol{\Omega(n)}$ & \textbf{\emph{n}} & $\boldsymbol{\Omega(n)}$ & \textbf{\emph{n}} & $\boldsymbol{\Omega(n)}$ & \textbf{\emph{n}} & $\boldsymbol{\Omega(n)}$   \\\hline
12 & 32 & 28 & 61 & 52 & 91 & 76 & 90  \\ \hline
20 & 47 & 44 & 61 & 68 & 77 & 92 & 61  \\ \hline
\end{tabular}
\end{center}
\label{table:NumberSRings4p}
\end{table}

In the special case that $p=2^k+1$ is a Fermat prime, the above formula simplifies to $\Omega(4p) = 15k+17$. For safe primes, we always have $\Omega(4p)=61$. In \tableref{table:NumberSRings4p} we list all integers of the form $4p$ less than 100.  

The proof of \thmref{thm:main} and \thmref{thm:main2} (which proofs can be found in \secref{sec:semiprime} and \secref{sec:4p}, respectively) can be summarized as following. By the Fundamental Theorem of Schur Rings over Cyclic Groups (due to Leung and Man \cite{LeungII, LeungI}), all Schur rings over cyclic groups belong to one of four families, which we call the \emph{traditional Schur rings}: namely, the trivial Schur ring, automorphic Schur rings, direct products of Schur rings, and wedge products of Schur rings (see the next section for definitions). Because these four families often overlap, special care is taken to ensure that an exact count is made. The trivial case, as the name suggests, is easy to consider. The families of wedge and direct products are considered recursively. The automorphic Schur rings are in one-to-one correspondence with the subgroups of  $\Aut(\z_{n}) \cong \prod_{i=1}^k \Aut(\z_{p_i^{e_i}})\cong \prod_{i=1}^k \z_{(p_i-1)p_i^{e_i-1}}$, where $n=\prod_{i=1}^r p_i^{e_i}$ is the prime factorization of $n$. 

In \cite{Ziv14}, Ziv-Av enumerates all Schur rings over small finite groups up to order 63. In \cite{CountingII}, the second author enumerates all Schur rings over cyclic groups up to order 400. In both cases, this was accomplished by computer software. In all instances, the two enumerations agree with the formulas found herein.

\section{Traditional Schur Rings}\label{sec:counting}

In this section, we remind the reader of important counting techniques introduced in \cite{CountingII}.

For a Schur ring $\S$ over $G$ and a subgroup $H\le G$, we say $H$ is a an $\S$-\emph{subgroup} if $H$ can be partitioned using the primitive sets of $\S$. Then $\S_H:=\S\cap \QH$ is a Schur ring over $\H$ and is called a \emph{Schur subring} of $\S$. We say a Schur ring $\S$ is \emph{primitive} if the only $\S$-subgroups are $1$ and $G$. 

For any subgroup $\H\le \Aut(G)$, let $G^\H$ denote the \emph{automorphic Schur ring} associated to $\H$. In the case $\H=1$, $G^1 = \QG$, the group algebra itself. For simplicity of notation, this Schur ring, called the \emph{discrete Schur ring}, is simply denoted as $G$ as the coefficient ring will provide little consequence. In the case that $G$ is abelian and $\H=\langle *\rangle$, we denote $G^\H$ as $G^\pm$. 

For any $n$, there is exactly one \emph{trivial Schur ring} over $G$, namely $G^0:=\Span\{1, G\smallsetminus 1\}$. The only primitive Schur ring over $\z_n$ when $n$ is composite is this trivial ring.  Thus, the trivial family is disjoint from the other three traditional families and contains exactly one ring. On the other hand, all Schur rings over $\z_p$, where $p$ is prime, are primitive. As observed in \cite{CountingII, Counting}, all Schur rings over $\z_p$ are automorphic $\left(\z_p^0=\z_p^{\Aut(\z_p)}\right)$ and $\Omega(p)=x$, where $x$ is the number of divisors of $p-1$.

If $G=H\times K$, $\S$ is a Schur ring over $H$, and $\T$ is a Schur ring over $K$, then $\S\times \T=\S\otimes_\Q\T$ denotes the \emph{direct product} of $\S$ and $\T$. Note that both $H$ and $K$ are necessarily subgroups of $\S\times \T$. In fact, $(\S\times\T)_H=\S$ and $(\S\times \T)_K= \T$.  Let $H$ and $K$ be groups and let $\H\le \Aut(H)$ and $\K\le \Aut(K)$. If $G=H\times K$, then we may naturally view $\H\times \K$ as a subgroup of $\Aut(G)$. As observed in \cite{CountingII}, if $G=H\times K$ then $\S=\S_H\times \S_K$ is automorphic if and only if $\S_H$ and $\S_K$ are automorphic. 

A \emph{section} $U$ of $G$ is a pair of subgroups $U=[K,H]$ such that $1\le K\le H\le G$ and $K\trianglelefteq G$. We say that a section $U=[K,H]$ is \emph{proper} if $1<K\le H<G$. We say that a section is \emph{trivial} if $K=H$. As all subgroups of $\z_n$ are normal and uniquely determined by their orders, we shall denote the section $[\z_d, \z_e]$ simply as $[d,e]$. 

Given any proper section $U=[K,H]$, a Schur ring $\S$ over $H$, and a Schur ring $\T$ over $K$, we form the \emph{wedge product} $\S\wedge_U \T$ by constructing the common refinement of $\D(\S)$ and $\D(\pi^{-1}(\T))$, where $\pi : G\to G/K$ is the natural map. To guarantee that $(\S\wedge_U\T)_H=\S$ and $\pi(\S\wedge_U\T) = \T$, the extra compatibility condition $\pi(\S)= \T_{H/K}$ is required in this construction. When $U=[H,H]$ is trivial, compatibility is automatic. We say a Schur ring is \emph{wedge-decomposable} when there exists a proper section such that the Schur ring can be expressed as a wedge product of two other Schur rings. In particular, $H$ and $K$ are $\S\wedge_U\T$-subgroups. Otherwise, we say the Schur ring is \emph{wedge-indecomposable}. 

As explained in \cite{CountingII}, every Schur ring has a \emph{wedge core}, which is a maximal wedge-indecomposable Schur subring. Let $\Omega(n,\S)$ be the number of Schur rings over $\z_n$ for which $\S$ is its wedge core. Let $\Omega(n,d) = \Omega(n,\z_d)$. Clearly, $\Omega(n) = \sum_\S \Omega(n,\S)$, where the sum ranges over the indecomposable Schur subrings of $\z_n$. We note that the indecomposable Schur rings are necessarily trivial, direct products of indecomposable Schur rings, or automorphic (of course, not every automorphic ring is indecomposable). The first case is trivial to count and the second is also enumerated recursively.

Lastly, to count automorphic Schur rings over $\z_n$, we note that there is a one-to-one correspondence between subgroups of $\Aut(\z_n)$ and automorphic Schur rings. So, it suffices to count the number of subgroups of $\Aut(\z_n)$, that is, compute $|\L(\Aut(\z_n))|$. The problem of counting the number of subgroups of an abelian group is a well studied problem in the literature (for example, see the references in \cite{CountingII}). The following lemma due to C\u{a}lug\u{a}reanu \cite{Calug} will be useful in our enumeration of automorphic Schur rings. 


\begin{Lem}\label{lem:latticeabelian} The number of subgroups of $\z_{p^k}\times \z_{p^\ell}$ is given as
\[\left|\L\middle(\z_{p^k}\times \z_{p^\ell}\middle)\right| = \sum_{j=0}^{\min(k,\ell)} \phi(p^j)(k-j+1)(\ell-j+1),\] where $\phi$ denotes Euler's totient function.
\end{Lem}

\section{Proof of \thmref{thm:main}}\label{sec:semiprime}
We proceed to the proof of \thmref{thm:main}. Consider $\z_{pq}$, where $p=\prod_{i=1}^n r_i^{k_i} +1$ and $q=\prod_{i=1}^n r_i^{\ell_i} +1$ are primes. By the Fundamental Theorem, all the Schur rings over this group are traditional. So, we now consider each of these families. As observed above, since $pq$ is composite, the trivial Schur ring contributes one to the count of $\Omega(pq)$ and does not intersect the other three families.  Additionally, as all Schur rings over $\z_p$ and $\z_q$ are necessarily automorphic, all direct product rings over $\z_{pq}$ are automorphic and will be counted under the automorphic family.

We consider next the wedge-decomposable Schur rings over $\z_{pq}$. Over $\z_{pq}$, the only possible proper sections are $[p, p]$ and $[q, q]$, which are trivial. In the first case, any possible Schur ring over $\z_p$ could be wedged with any possible Schur ring over $\z_q$. This produces $\Omega(p)\Omega(q)=xy$ distinct Schur rings. The subgroups of each of these Schur rings will be exactly $1$, $\z_p$, and $\z_{pq}$. Note that $\z_q$ is missing since elements of order $q$ and $pq$ are fused together in cosets of $\z_p$. Thus, these Schur rings are not automorphic. The second case is similar and produces $xy$ distinct Schur rings from the automorphic ring and the wedge products already accounted for. Therefore, there are $2xy$ wedge products over $\z_{pq}$. 

Note that by the above decomposition of $p$, $p-1$ has $\prod_{i=1}^n(k_i+1)$ divisors. Likewise, $q-1$ has $\prod_{i=1}^n(\ell_1+1)$ divisors. Considering the wedge-decomposable and trivial families, we have already accounted for $2\prod_{i=1}^n (k_i+1)(\ell_i+1)+1$ distinct, non-automorphic Schur rings. Thus, to prove \thmref{thm:main} it suffices to count the number of distinct automorphic Schur rings over $\z_{pq}$. As these Schur rings are in direct one-to-one correspondence with the subgroups of $\Aut(\z_{pq})$, we see that 
\begin{equation}\label{eq:mainalmost}\Omega(pq) = |\L(\Aut(\z_{pq}))| + 2\prod_{i=1}^n (k_i+1)(\ell_i+1)+1.\end{equation} 

For a finite group $G$, let $G=\prod_{i=1}^k P_i$ be its primary decomposition. Then it is well-known that $\L(G) \cong \prod_{i=1}^k \L(P_i)$. In the case of $\z_{pq}$, we see that $\Aut(\z_{pq}) \cong \z_{p-1}\times \z_{q-1}$. Hence, the primary decomposition for $\Aut(\z_{pq})$ is given as
\[\Aut(\z_{pq}) \cong \prod_{i=1}^n (\z_{r_i^{k_i}}\times \z_{r_i^{\ell_i}}).\] By \lemref{lem:latticeabelian},  \[\left|\L\middle(\z_{r_i^{k_i}}\times \z_{r_i^{\ell_i}}\middle)\right| = \sum_{j=0}^{\min(k_i,\ell_i)}\phi(r_i^j)(k_i-j+1)(\ell_i-j+1).\] Therefore, \begin{equation}\label{eq:mainmissinglink}|\L(\Aut(\z_{pq}))|=\prod_{i=1}^n\sum_{j=0}^{\min(k_i,\ell_i)} \phi(r_i^j)(k_i-j+1)(\ell_i-j+1).\end{equation} Finally, \thmref{thm:main} follows immediately from \eqref{eq:mainalmost} and \eqref{eq:mainmissinglink}, which finishes the proof.

\begin{Exam} We present a complete enumeration of the Schur rings over $\z_{21}$ as an example to illustrate the previous proof. There are $\Omega(3)=2$ Schur rings over $\z_3$, namely $\z_3^0$ and $\z_3$. There are $\Omega(7)=4$ Schur rings over $\z_7$, namely, $\z_7^0$, $\z_7^{\langle 2\rangle}$, $\z_7^\pm$, and $\z_7$. 

Below we list the $\Omega(21)=27$ Schur rings over $\z_{21}$:
\[\z_3^0 \wedge \z_7^0, \z_3^0 \wedge \z_7^{\langle 2\rangle}, \z_3^0 \wedge \z_7^\pm, \z_3^0 \wedge \z_7, \z_3 \wedge \z_7^0, \z_3 \wedge \z_7^{\langle 2\rangle}, \z_3 \wedge \z_7^\pm, \z_3 \wedge \z_7, \]
\[\z_7^0 \wedge \z_3^0, \z_7^{\langle 2\rangle} \wedge \z_3^0, \z_7^\pm \wedge \z_3^0, \z_7 \wedge \z_3^0, \z_7^0 \wedge \z_3, \z_7^{\langle 2\rangle} \wedge \z_3, \z_7^\pm \wedge \z_3, \z_7 \wedge \z_3, \]
\[\z_3^0 \times \z_7^0, \z_3^0 \times \z_7^{\langle 2\rangle}, \z_3^0 \times \z_7^\pm, \z_3^0 \times \z_7, \z_3 \times \z_7^0, \z_3 \times \z_7^{\langle 2\rangle}, \z_3 \times \z_7^\pm, \z_3 \times \z_7\ (\cong \z_{21}), \]
\[\z_{21}^0, \z_{21}^{\langle 5 \rangle},\z_{21}^\pm. \qedhere\]
\end{Exam}

\section{Proof of \thmref{thm:main2}}\label{sec:4p}




Let $p=2^ka+1$, where $a$ is odd. Let $x$ denote the number of divisors of $p-1$. Then $x/(k+1)$ is the number of divisors of $a$. 

The construction of wedge products will be more complicated in the case $n=4p$ as wedge-compatibility is not necessarily trivial. There is also the concern that different sections can potentially create the same wedge product, e.g. if $U=[2,4]$ and $U'=[4,4]$, then $\z_4 \wedge_U \z_6 = \z_4\wedge_{U'} \z_3$. Additionally, the wedge product construction is associative, e.g. $(\z_2\wedge\z_2)\wedge \z_3 = \z_2\wedge(\z_2\wedge \z_3)$. 
To avoid these issues, we will restrict our attention to wedge-decomposable Schur rings $\S\wedge_U \T$ where $U=[K,H]$, $K$ is a minimal $\S$-subgroup, and $\S$ is wedge-indecomposable, that is, $\S$ is the wedge-core of $\S\wedge_U \T$. 

We will organize Schur rings over $\z_{4p}$ according to the order of its wedge-core.  We need to consider $\Omega(4p,d)$ for the divisors $d=2,p,4,2p$ ($d=1, 4p$ are omitted as these are not proper sections). If $\S$ is primitive, the only possible section would be trivial. As such, $\Omega(n,\S)=\Omega(n/d)$. For example, $\Omega(4p,2)=\Omega(2p)= 3x+1$ and  $\Omega(4p, p)=\Omega(4)=3$. In fact, as all the $x$ many Schur rings $\S$ over $\z_p$ are primitive, $\Omega(4p,\S)=\Omega(4)$. Hence, there are $\Omega(p)\Omega(4) = 3x$ many Schur rings over $\z_{4p}$ with a core of order $p$.

The three Schur rings over $\z_4$ are $\z_4^0$, $\z_2\wedge \z_2$, and $\z_4$, where the first and last ones are indecomposable. 
As $\z_4^0$ is primitive, $\Omega(4p, \z_4^0) = \Omega(p) = x$. For $\S=\z_4$ we select the section $U=[2,4]$. Hence, $\Omega(4p, 4)$ counts the number of Schur rings of the form $\z_4\wedge_U \T$ where $\T$ is a Schur ring over $\z_{2p}$ and $\z_2$ is a subring. To count the number of such rings $\T$ we revisit the proof of \thmref{thm:main} with the simplification $q=2$, which shows that there are exactly $x$ many Schur rings of the form $\z_2\wedge \S$, $x$ many Schur rings of the form $\S\wedge\z_2$, $x$ many Schur rings of the form $\z_2\times \S$, and one trivial ring $\z_{2p}^0$. Hence, there are $2x$ many Schur rings over $\z_{2p}$ that contain the subring $\z_2$, that is, $\Omega(4p, 4) = 2x$. Therefore, there are $3x$ many Schur rings over $\z_{4p}$ with a core of order 4. 

It is also observed that there are $x+1$ indecomposable Schur rings over $\z_{2p}$, which $x$ have the form $\z_2\times \S'$ and one is trivial. Note $\Omega(4p, \z_{2p}^0)=\Omega(2)=1$, that is, $\z_{2p}^0\wedge \z_2$ is the only option.  If  $\S=\z_2\times \S'$ is the core, then $\S$ has two distinct minimal subgroups, namely, $\z_2$ and $\z_p$. As such, two proper sections need to be considered, $[2,2p]$ and $[p,2p]$. For fixed $\S'$, there are $2$ Schur rings over $\z_{2p}$ which contain $\S'$ as a subring, namely $\S'\wedge \z_2$ and $\S'\times \z_2$. Hence, there are two possibilities for $(\z_2\times \S')\wedge_{[2,2p]} \T$. For $(\z_2\times \S')\wedge_{[p,2p]} \T$, we count the number of Schur rings over $\z_4$ which contain $\z_2$ as a subgroup. There are two such rings, $\z_2\wedge\z_2$ and $\z_4$. Lastly, there is  $\Omega(2)=1$ ring of the form $(\z_2\times \S')\wedge_{[2p,2p]} \T$, namely $(\z_2\times \S')\wedge \z_2$. Therefore, there are $2+2-1=3$ many Schur rings over $\z_{4p}$ with $\z_2\times \S'$ as its core. This accounts for $3x+1$ many Schur rings over $\z_{4p}$ with a core of order $2p$.

In summary, this accounts for $(3x+1)+3x+3x+(3x+1) = 12x+2$ many wedge-decomposable Schur rings over $\z_{4p}$. It remains to consider the indecomposable Schur rings over $\z_{4p}$. There is, of course, the trivial Schur ring $\z_{4p}^0$, as well as the indecomposable automorphic Schur rings and those direct products of the form $\z_4^0\times \S$ for some Schur rings $\S$ of order $p$. In regard to the automorphic Schur rings, we know the total count is equal to $|\L(\Aut(\z_{4p}))|$. The sublattice $\L(\Aut(\z_4))\times \L(\Aut(\z_{p}))$ will consist of two lattice-isomorphic copies of $\L(\z_{p-1})$. The full lattice $\L(\Aut(\z_{4p}))$ contains these two layers and all the diagonal entries that sit between the top and bottom layers. Those non-diagonal automorphic Schur rings in the top layer have the form $\z_4\times \S$, for some Schur ring $\S$ of order $p$, and are indecomposable.  Those non-diagonal automorphic Schur rings in the bottom layer have the form $(\z_2\wedge\z_2)\times \S$ and are wedge-decomposable. Notice that these decomposable automorphic Schur rings are in one-to-one correspondence with those indecomposable Schur rings of the form $\z_4^0\times \S$. Thus, if every diagonal automorphic Schur ring over $\z_{4p}$ is indecomposable, which we claim, then the number of indecomposable Schur rings over $\z_{4p}$ is $1+|\L(\Aut(\z_{4p}))|$. Using \lemref{lem:latticeabelian}, we see that $\Aut(\z_{4p})\cong \z_2 \times \z_{p-1} \cong (\z_2\times \z_{2^k})\times \z_a$ and 
\[|\L(\Aut(\z_{4p}))|  = |\L(\z_2\times \z_{2^k})||\L(\z_a)| = (2(k+1)+k)\left(\dfrac{x}{k+1}\right)=\dfrac{3k+2}{k+1}x.\]

To prove the claim, we introduce the representation $\omega : \Q[\z_n] \to \Q(\zeta_n)$ which maps a generator of $\z_n$ to $\zeta_n:=e^{2\pi i/n}$. We remind the reader that in \cite{CountingII} we saw that an automorphic Schur ring $\S$ is wedge-decomposable if and only if $\omega(\S)\le \Q(\zeta_d)$ for some proper divisor $d\mid n$ (excluding, of course, the case when $n$ is prime). By definition, the diagonal automorphic Schur rings are not contained in $\L(\Aut(\z_4))\times \L(\Aut(\z_p))$, which implies that their image is not contained in $\L(\Q(\zeta_4))$ or  $\L(\Q(\zeta_p))$. This implies they are all indecomposable, as claimed.

As we have now exhausted all possibilities, we see that 
\[\Omega(4p) = (12x+2) + 1+\dfrac{3k+2}{k+1}x = \dfrac{15k+14}{k+1}x+3,\] 
which finishes the proof of \thmref{thm:main2}.

\begin{Exam} We present a complete enumeration of the Schur rings over $\z_{12}$ as an example to illustrate the previous proof. Note $12=4(3)$ and $3=2^1\cdot 1 + 1$. There are $\Omega(3)=2$ Schur rings over $\z_3$, namely $\z_3^0$ and $\z_3$. There are $\Omega(6)=7$ Schur rings over $\z_{6}$, namely 
\[\z_2\wedge \z_3^0, \z_2\wedge \z_3, \z_3^0\wedge \z_2, \z_3\wedge \z_2, \z_6^0, \z_2\times \z_3^0\ (=\z_6^\pm),  \z_2\times \z_3\ (=\z_6). \]

Below we list the $\Omega(12)=32$ Schur rings over $\z_{12}$:
\[\z_2\wedge \z_2\wedge \z_3^0, \z_2\wedge \z_2\wedge \z_3, \z_2\wedge \z_3^0\wedge \z_2, \z_2\wedge \z_3\wedge \z_2, \z_2\wedge \z_6^0, \z_2\wedge \z_6^\pm, \z_2\wedge \z_6, \]
\[\z_3^0 \wedge \z_4^0, \z_3^0 \wedge \z_2\wedge \z_2, \z_3^0 \wedge \z_4, \z_3 \wedge \z_4^0, \z_3 \wedge \z_2\wedge \z_2, \z_3 \wedge \z_4, \]
\[\z_4^0 \wedge \z_3^0, \z_4^0 \wedge \z_3, \z_4 \wedge_{[2,4]} (\z_2\wedge \z_3^0), \z_4\wedge_{[2,4]} (\z_2\wedge \z_3),  \z_4\wedge_{[2,4]} \z_6^\pm, \z_4\wedge_{[2,4]} \z_6, \]
\[\z_6^0 \wedge \z_2, \z_6^\pm \wedge \z_2, \z_6^\pm \wedge_{[2,6]} \z_6^\pm, \z_6^\pm \wedge_{[3,6]} \z_4, \z_6 \wedge \z_2, \z_6 \wedge_{[2,6]} \z_6,  \z_6 \wedge_{[3,6]} \z_4, \]
\[\z_{12}^0, \z_4^0\times \z_3^0, \z_4^0\times \z_3, \z_4\times \z_3^0, \z_4\times \z_3\ (=\z_{12}), \z_{12}^\pm. \qedhere\]

\end{Exam}

\bibliographystyle{plain}
\bibliography{Srings}

\begin{thebibliography}{1}

\bibitem{Calug}
Grigore {C\u{a}lug\u{a}reanu}.
\newblock {The total number of subgroups of a finite Abelian group.}
\newblock {\em {Sci. Math. Jpn.}}, 60(1):157--167, 2004.

\bibitem{LeungII}
Ka~Hin {Leung} and Shin~Hing {Man}.
\newblock On {S}chur rings over cyclic groups {II}.
\newblock {\em {Journal of Algebra}}, 183:273--285, 1996.

\bibitem{LeungI}
Ka~Hin {Leung} and Shin~Hing {Man}.
\newblock On {S}chur rings over cyclic groups.
\newblock {\em {Israel Journal of Mathematics}}, 106:251--267, 1998.

\bibitem{CountingII}
Andrew {Misseldine}.
\newblock On counting {S}chur rings over cyclic groups.
\newblock pre-print (submitted).

\bibitem{MePhD}
Andrew {Misseldine}.
\newblock {\em Algebraic and Combinatorial Properties of {S}chur Rings over
  Cyclic Groups}.
\newblock PhD thesis, Brigham Young University, 2014.

\bibitem{Counting}
Andrew {Misseldine}.
\newblock Counting {S}chur rings over cyclic groups.
\newblock {\em {Journal of Algebraic Combinatorics}}, 51:155--169, Feb 2020.

\bibitem{Muzychuk09}
Mikhail {Muzychuk} and Ilia {Ponomarenko}.
\newblock {S}chur rings.
\newblock {\em {European Journal of Combinatorics}}, 30:1526--1539, 2009.

\bibitem{Ziv14}
Matan {Ziv-Av}.
\newblock Enumeration of schur rings over small groups.
\newblock In {\em {Computer algebra in scientific computing. 16th international
  workshop, CASC 2014, Warsaw, Poland, September 8--12, 2014. Proceedings}},
  pages 491--500. Berlin: Springer, 09 2014.

\end{thebibliography}
\end{document}